\documentclass{amsart}

\numberwithin{equation}{section}

\begin{document}
\title[]{Response to a setpoint change in PID controlled time-delay feedback systems}


\author{Gianpasquale Martelli}
\address{Via Domenico da Vespolate 8 28079 Vespolate Italy}
\curraddr{}
\email{gianpasqualemartelli@libero.it}
\thanks{}



\date{February 22, 2007}

\dedicatory{}

\begin{abstract}
The response to a setpoint change in PID controlled feedback systems plays an important role for the tuning methods. This response may be easily evaluated in linear systems without delay by solving the related ordinary differential equations. These equations are also obtained for system with delay, if the  delay exponential term is approximated to a polynomial by one of the   Pad\'e expressions. If this simplification is avoided, one must deal with differential difference equations. In this paper their explicit solution, obtained with the method of steps, is presented. Moreover a second order equation is investigated as example and the complete set of the coefficients formulas is given.
\end{abstract}

\maketitle

\section{Introduction}
The closed-loop response of a feedback control system can be considered one of the main behaviour taken into consideration for a proper tuning. A very common and normally desired goal consists of a deadbeat response (see para. 10.2.5 of \cite{bib1} p. 171) i. e. the wanted value of the controlled variable, after a setpoint change, is achieved as fast as possible and held within a narrow band. This operation mode has been completely studied and successfully applied only to the systems without delay.
The Ziegler-Nichols tuning method \cite{bib2}, which is the most known and used one, considers the quarter amplitude decay ratio, obtains the process parameters values from an open-loop test and consequently gives simple relationships between them and the  controller ones. An another is the Kappa-Tau method  \cite{bib3}, based on a step response experiment, which supplies parameters relationships obtained applying dominant pole design, at constant maximum sensivity gain, on typical process control models. Moreover there are optimization criteria, minimizing the integrated error during a step response, but calculated after the reduction of the exponential time-delay term to a polynomial form by means of a Pad\'e approximation (see for example \cite{bib4}). Finally the actual commercial controllers, which can be single in separate hardware box or multiple in distributed control systems (DCS), are provided with adaptive techniques (see para. 52.4 of \cite{bib1} p. 825). Most of these techniques are based on step response analysis and the controller parameters are adjusted according to some rules in order to have a performance index inside a preset range.

From these preliminary remarks it follows that a strict and simple formula, not yet available, of a setpoint change response can be helpful  for an accurate and safe controller tuning of time-delay systems and hence it will be the aim of this paper.

The transfer functions of the PID controller $T_{c}(s)$ and process $T_{p}(s)$ , taken into consideration, are given by
\begin{equation}\label{eq:1.1}
T_{c}(s)=k+k_{i}/s+k_{d} \, s
\end{equation}
\begin{equation}\label{eq:1.2}
T_{p}(s)=e^{-s} \frac{\sum_{h=0}^{h=m_{b}-2}e_{h} s^{h}} {\sum_{h=0}^{h=m_{a}-1}a_{h} s^{h}}\  
\end{equation}
where  $m_{b}<=m_{a}$.
The differential difference equation related to the closed-loop tranfer function is described by
\begin{equation}\label{eq:1.3}
\sum_{h=1}^{h=m_{a}}a_{h} \frac{d^{h}y(t)} {dt^{h}}\  = - \sum_{h=0}^{h=m_{b}}b_{h} \frac{d^{h}y(t-1)} {dt^{h}}\ + \sum_{h=0}^{h=m_{c}}c_{h} \frac{d^{h}f(t-1)} {dt^{h}}\
\end{equation}
subject to an initial condition of the form $y(t)=y_{0,k}(t)$ for each of the $q$ consecutive intervals  $-1+ \tau_{k-1} \leq t \leq -1+ \tau_{k}$ ($ \tau_{k-1}< \tau_{k}$, $ \tau_{0}=0$ and $ \tau_{q}=1$ ) and where
\begin{itemize}
\item[-] $y(t)$ is the controlled variable
\item[-] $f(t)$ is the setpoint
\item[-] t is the normalized time referred to the time delay
\item[-] $a_{h}$, $b_{h}$ and  $c_{h}=b_{h}$ are constants
\item[-] $m_{a}$, $m_{b}$ and $m_{c}=m_{b}$ are positive integers
\end{itemize}

The solution of (\ref{eq:1.3}) can be evaluated by means of the method of the Laplace transform, as a contour integral or as an infinite exponential series,  or by means of the method of steps \cite{bib5}, \cite{bib6}. The contour integral cannot be expressed in terms of elementary functions and for the exponential series, helpful for the study of the behaviour for higher values of the independent variable, an approximation is needed; therefore both are not suitable for an accurate evaluation of the solution for low values of the time. On the contrary this is not true for the method of steps, since, denoting $y_{n,k}$ as the solution valid for $n-1+ \tau_{k-1} \leq t \leq n-1+ \tau_{k}$, $y_{1,k}$ can be analytically evaluated solving (\ref{eq:1.3}), which can be considered as an ordinary differential equation. Any consecutive $y_{n,k}$  can be analogously determined, but tedious calculations will be necessary at each time if a general expression will not be found.

This general solution exists provided that the initial condition $y_{0,k}(t)$ and the forcing term $f_{n,k}(t)$ are each a finite sum of the terms of the complementary solution of (\ref{eq:1.3}), but with a null time delay, each multiplied by a polynomial of $t$.
The solution $y_{n,k}(t)$ will be of the same form. In detail:
\begin{equation}\label{eq:1.4}
y_{n,k}(t) = \sum_{p=1}^{p=m_{a}} e^{r_{p}t} \sum_{i=0}^{i=v_{k,p}+n} G^{*}_{n,k,p,i}t^{i}
\end{equation}
\begin{equation}\label{eq:1.5}
f_{n,k}(t) =  \sum_{p=1}^{p=m_{a}} e^{r_{p}t} \sum_{i=0}^{i=w_{k,p}+n} F^{*}_{n,k,p,i}t^{i}
\end{equation}
where $v_{k,p}$ and $w_{k,p}$ are positive integers with $v_{k,p}>=w_{k,p}$, $1<=k<=q$ and  $r_{p}$ are the characteristic roots of the corresponding ordinary differential equation without time delay. These roots are assumed all simple and real in this paper, but the adopted procedure can be used also if multiple real or complex roots are present.  

It is worthwhile to note that, thanks to the integrating term of the PID controller, one of the characteristic roots $r_{p}$ is null and therefore steady initial conditions and setpoint step changes, which are usually encountered in the tuning methods, yield solutions compatible with (\ref{eq:1.4}).

Since  the coefficients  $G^{*}_{u,k,p,i}$ are considered, in the method of steps, unknown for $u=n$ and known for  $u=n-1$, the necessary and sufficient conditions for the correctness of (\ref{eq:1.4}), for a given $u$, are:
\begin{enumerate}
\item equal number of each different term depending on the time $e^{r_{p}t}t^{i}$ included in both sides of  (\ref{eq:1.1}), modified by the introduction of (\ref{eq:1.4})  
\item the total number of the coefficients  $G^{*}_{u,k,p,i}$  equal to the sum of the number of all above mentioned terms $e^{r_{p}t}t^{i}$ and the continuity equations one.
\end{enumerate}
These two conditions are fulfilled and will be checked in the next Section.

\section{Main results}
The sought expressions of the coefficients are simpler if, for each $y_{n,k}$ and  $f_{n,k}$, the origin of the temporal axis $t_{n}$ is assumed coincident with $n-1$ ($t = n-1+t_{n}$). Therefore (\ref{eq:1.4}) and (\ref{eq:1.5}) can be rewritten, for $n>=0$, as
\begin{equation}\label{eq:2.1}
y_{n,k}(t_{n}) = \sum_{p=1}^{p=m_{a}} e^{r_{p}t_{n}} \sum_{i=0}^{i=v_{k,p}+n} G_{n,k,p,i}t_{n}^{i}
\end{equation}
\begin{equation}\label{eq:2.2}
f_{n,k}(t_{n}) = \sum_{p=1}^{p=m_{a}} e^{r_{p}t_{n}} \sum_{i=0}^{i=w_{k,p}+n} F_{n,k,p,i}t_{n}^{i}
\end{equation}

Since $t_{n-1}-1=t_{n}$ holds and since the trasformations from (\ref{eq:1.4}) to (\ref{eq:2.1}) and from (\ref{eq:1.5}) to (\ref{eq:2.2}) substitute $t$ with $t_{n}$ in $y_{n,k}$ and with $t_{n-1}$ in $y_{n-1,k}$ and $f_{n-1,k}$, (\ref{eq:1.3}) becomes
\begin{equation}\label{eq:2.4}
\sum_{h=1}^{h=m_{a}}a_{h} \frac{d^{h}y_{n,k}(t_{n})} {dt_{n}^{h}}\  = - \sum_{h=0}^{h=m_{b}}b_{h} \frac{d^{h}y_{n-1,k}(t_{n})} {dt_{n}^{h}}\ + \sum_{h=0}^{h=m_{c}}c_{h} \frac{d^{h}f_{n-1,k}(t_{n})} {dt_{n}^{h}}\
\end{equation}
subject to the following continuity equations:
\begin{equation}\label{eq:2.5}
\begin{split}
y_{n,1}(0)= &y_{n-1,q}(1) \quad for \,k=1 \\
y_{n,k}( \tau_{k-1})= &y_{n,k-1}( \tau_{k-1}) \quad for \,k>1
\end{split}
\end{equation}
\begin{equation}\label{eq:2.6}
\begin{split}
\frac{d^{h}y_{n,1}(0)} {dt_{n}^{h}}\ = &\frac{dy^{h}_{n-1,q}(1)}{dt_{n}^{h}}\ \quad for \,k=1\\ \frac{d^{h}y_{n,k}( \tau_{k-1})} {dt_{n}^{h}}\ =& \frac{d^{h}y_{n,k-1}(\tau_{k-1})}{dt_{n}^{h}}\  \quad for \,k>1 \, and \, 0<h<m_{a}
\end{split}
\end{equation}

Let introduce in  (\ref{eq:2.4}) the solution (\ref{eq:2.1}) and  its derivatives, evaluated according to (\ref{eq:A.1}). The necessary and sufficient conditions of the correctness of (\ref{eq:2.1})  can be detailed for each $n>0$  and each $k$ as follows:

\begin{enumerate}
\item number of the terms $e^{r_{p}t}t^{i}$
\begin{enumerate}
\item left-hand side of (\ref{eq:2.4}): $ \sum_{p=1}^{p=m_{a}}(v_{k,p}+n)$\\
 it is $0 \leq i \leq v_{k,p}+n$, but the coefficient of $t_{n}^{v_{k,p}+n}$ is null since it is equal to the characteristic  polynomial of the corresponding equation without time delay
\item right-hand side of (\ref{eq:2.4}): $ \sum_{p=1}^{p=m_{a}}(v_{k,p}+n)$\\
 it is $0 \leq i \leq v_{k,p}+n-1$
\end{enumerate}
\item  coefficients $G_{n,k,p,i}$  ( $1 \leq i \leq v_{k,p}+n$)
\begin{enumerate}
\item number of the equations (\ref{eq:2.4}) (see item (1)): $ \sum_{p=1}^{p=m_{a}}(v_{k,p}+n)$
\item number of the unknowns : $ \sum_{p=1}^{p=m_{a}}(v_{k,p}+n)$ \\
$G_{n,k,p,0}$ has as coefficient the characteristic polynomial of the corresponding equation without time delay and then disappears
\end{enumerate}
\item  coefficients $G_{n,k,p,0}$ ( $1 \leq p \leq m_{a}$)
\begin{enumerate}
\item number of the equations (\ref{eq:2.5}) and (\ref{eq:2.6}): $m_{a}$
\item number of the unknowns: $m_{a}$
\end{enumerate}

\end{enumerate}

The analytical trasformations required for the evaluation of the coefficients formulas are quite simple. Since  series products exist in both sides of (\ref{eq:2.4}) and the  series of the power of $t$ is in the last right position, it is only necessary to shift this series from its actual position to the first left one, modifying accordingly the initial and the final values of all indexes. In this manner the coefficients of $e^{r_{p}t}t^{i}$ of both sides, each consisting of the entire expression following this first series, can be correctly equated and therefore  (\ref{eq:2.4}) become an identity for any value of $t$. 

In detail:
\begin{enumerate}
\item $i>0$\\
Introducing (\ref{eq:2.1}) and  its derivatives, evaluated according to (\ref{eq:A.1}), in the left-hand side $LHS$ and in the  right-hand side $RHS$ of (\ref{eq:2.4}),  we obtain
\begin{equation}\label{eq:2.7}
LHS =  \sum_{p=1}^{p=m_{a}} (S_{1}(m_{a)}+S_{2}(m_{a},v_{k,p}+n)) P_{n,k,h,p,i,j} e^{r_{p}t_{n}}  t_{n}^{j}
\end{equation}
\begin{equation}\label{eq:2.8}
\begin{split}
RHS = &+ \sum_{p=1}^{p=m_{a}} \sum_{j=0}^{j=v_{k,p}+n-1} (Q_{n,k,0,p,j,j}+  R_{n,k,0,p,j,j}) e^{r_{p}t_{n}}  t_{n}^{j}\\ &+ \sum_{p=1}^{p=m_{a}}  (S_{1}(m_{b)}+S_{2}(m_{b},v_{k,p}+n-1))   Q_{n,k,h,p,i,j} e^{r_{p}t_{n}}  t_{n}^{j}\\ &+ \sum_{p=1}^{p=m_{a}}   (S_{1}(m_{c)}+S_{2}(m_{c},w_{k,p}+n-1)) R_{n,k,h,p,i,j} e^{r_{p}t_{n}}  t_{n}^{j}
\end{split}
\end{equation}
where
\begin{displaymath}
S_{1}(m) = \sum_{h=1}^{h=m} \sum_{i=0}^{i=h-1} \sum_{j=0}^{j=i}
\end{displaymath}
\begin{displaymath}
S_{2}(m,z) = \sum_{h=1}^{h=m} \sum_{i=h}^{i=z} \sum_{j=i-h}^{j=i}
\end{displaymath}
\begin{displaymath}
P_{n,k,h,p,i,j} =a_{h} r_{p}^{h-i+j} \frac{i!h!}{j!(i-j)!(h-i+j)}\ G_{n,k,p,i}
\end{displaymath}
\begin{displaymath}
Q_{n,k,h,p,i,j} =- b_{h} r_{p}^{h-i+j} \frac{i!h!}{j!(i-j)!(h-i+j)}\ G_{n-1,k,p,i} 
\end{displaymath}
\begin{displaymath}
R_{n,k,h,p,i,j} =c_{h} r_{p}^{h-i+j} \frac{i!h!}{j!(i-j)!(h-i+j)}\ F_{n-1,k,p,i}
\end{displaymath}
It is now necessary to shift the index $j$ series, in $S_{1}(m)$ and $S_{2}(m,z)$, from the last right position to the first left one according to the series identities detailed in Appendix B. Before making equal the coefficients of  $t^{j}$ included in $LHS$  (\ref{eq:2.7}) and in $RHS$  (\ref{eq:2.8}) it is worthwhile to note that, for $i=j$, $P_{n,k,h,p,i,j}$ is given by
\begin{displaymath}
P_{n,k,h,p,i,i} =a_{h} r_{p}^{h} G_{n,k,p,i} 
\end{displaymath}
and that the following holds
\begin{equation}\label{eq:2.9}
 \sum_{h=1}^{h=m_{a}} P_{n,k,h,p,i,i} = G_{n,k,p,i} \sum_{h=1}^{h=m_{a}} a_{h} r_{p}^{h} = 0
\end{equation}
because the coefficient of $G_{n,k,p,i}$, which is the characteristic polynomial of the corresponding ordinary differential equation without time delay, is null. 
The coefficients of the terms $t^{j}$ of $LHS$ (\ref{eq:2.7}) are consequently  deleted as follows:
\begin{enumerate}
\item term with $j=0$:  $G_{n,k,p,0}$, existing only in this case
\item each term from $j=1$ to $j=n+v_{k,p}-1$: $G_{n,k,p,j}$ 
\item term with  $j=n+v_{k,p}$: $G_{n,k,p,n+v_{k,p}}$, only one present 
\end{enumerate}
Finally the coefficient matrix of the system, consisting of $v_{k,p}+n$ linear equations obtained making equal the coefficients of the terms $t^{j}$ belonging to  $LHS$  and  $RHS$, is triangular and provided with  $m_{a}$ diagonals starting from the central one. The unknowns $G_{n,k,p,i}$, from  $i=n+v_{k,p}$ to  $i=1$,  can be evaluated in function of $G_{n-1,k,p,i}$, for each value of $n$, $k$ and $p$, recursively or with the analytical solution of the system, obtained by means of a consecutive elimination of the unknowns.

\item  $i=0$\\ Defining $H_{e,n,k,h}(t_{n})$ as 
\begin{displaymath}
\begin{split}
H_{e,n,k,0}(t_{n}) =  \sum_{p=1}^{p=m_{a}}  e^{r_{p}t_{n}} \sum_{i=e}^{i=n+v_{k,p}}  G_{n,k,p,i} t_{n}^{i} \quad for \, h=0
\end{split}
\end{displaymath}
\begin{displaymath}
\begin{split}
H_{e,n,k,h}(t_{n}) = &\frac{d^{h}} {dt^{h}}\ \left(   \sum_{p=1}^{p=m_{a}}  e^{r_{p}t_{n}} \sum_{i=e}^{i=n+v_{k,p}}   G_{n,k,p,i} t_{n}^{i} \right)\\ = &\sum_{p=1}^{p=m_{a}} \left(   \sum_{i=e}^{i=h-1}   \sum_{j=0}^{j=i} + \sum_{i=h}^{i=n+v_{k,p}}   \sum_{j=i-h}^{j=i} \right)  r_{p}^{h-i+j}\\ & \frac{i!h!}{j!(i-j)!(h-i+j)}\ G_{n,k,p,i} e^{r_{p}t_{n}} t_{n}^{j} \quad \, for \, h>0
\end{split}
\end{displaymath}
 (\ref{eq:2.5}) and (\ref{eq:2.6}) become, for $ 0<=h<=m_{a}-1$, 
\begin{equation}\label{eq:2.10}
\begin{split}
\sum_{p=1}^{p=m_{a}}r_{p}^{h} G_{n,1,p,0} = &-H_{1,n,1,h}(0) +H_{0,n-1,q,h}(1)  \quad for \,k=1\\ \sum_{p=1}^{p=m_{a}} e^{r_{p} \tau_{k-1}} r_{p}^{h} G_{n,k,p,0} = &-H_{1,n,k,h}( \tau_{k-1})+H_{0,n,k-1,h}( \tau_{k-1})  \quad for \,k>1
\end{split}
\end{equation}
For each value of $k$ the coefficient matrix of the system, consisting of the $m_{a}$ linear equations  (\ref{eq:2.10}), is a Vandermonde matrix, whose terms of each column, for $k>1$, are multiplied by a constant. This system  can be analytically solved with respect to the $m_{a}$ unknowns $G_{n,k,p,0}$ by means of the Cramer's method and of the Laplace's expansion, thanks to the known expression of the Vandermonde determinant,  as detailed in Appendix C. 
\end{enumerate}

\section{Examples} 
Let consider a first-order plant and hence assume $m_{a}=2$, $m_{b}=2$, $a_{0}=0$, $a_{1}=1$, $a_{2}=t_{p}$. Two cases, whose the first related to a preset initial condition with $q>1$ and to a setpoint always null and the second related to  a steady initial condition with $q=1$ and to a setpoint change from $1$ to $0$, are examined.

\subsection{Example no.1} 
Since $q>1$, $f_{n,k}=0$ ($n>=0$),  (\ref{eq:2.4}) and  (\ref{eq:2.1}) become
\begin{equation}\label{eq:3.1}
\begin{split}
\frac{dy_{n,k}(t_{n})} {dt_{n}}\ &+t_{p} \frac{d^{2}y_{n,k}(t_{n})}{dt_{n}^{2}}\ =- b_{0} y_{n-1,k}(t_{n})\\ &-b_{1} \frac{dy_{n-1,k}(t_{n})}{dt_{n}}\ -b_{2} \frac{d^{2}y_{n-1,k}(t_{n})} {dt_{n}^{2}}\
\end{split}
\end{equation}
\begin{equation}\label{eq:3.2}
y_{n,k}(t_{n}) =  \sum_{i=0}^{i=v_{k,1}+n} G_{n,k,1,i}t_{n}^{i}+ e^{-t_{n}/t_{p}} \sum_{i=0}^{i=v_{k,2}+n} G_{n,k,2,i}t_{n}^{i}
\end{equation}

The expressions of $G_{n,k,p,i}$ for $n>0$ are given by
\begin{itemize}
\item[a)] $i>0$
\begin{displaymath}
G_{n,k,1,i}  =  - \frac {b_{0}}{i}\ G_{n-1,k,1,i-1}  \quad for \, i=v_{k,1}+n
\end{displaymath} 
\begin{displaymath}
\begin{split} 
G_{n,k,1,i} + t_{p}(i+1)G_{n,k,1,i+1} =  &- \frac {b_{0}}{i}\
G_{n-1,k,1,i-1}-b_{1}G_{n-1,k,1,i}\\ &for \, i=v_{k,1}+n-1
\end{split}
\end{displaymath} 
\begin{displaymath}
\begin{split} 
G_{n,k,1,i} &+ t_{p}(i+1)G_{n,k,1,i+1} =  - \frac {b_{0}}{i}\ G_{n-1,k,1,i-1}-b_{1}G_{n-1,k,1,i} \\  &- b_{2}(i+1)G_{n-1,k,1,i+1} \quad for \, i<=v_{k,1}+n-2
\end{split}
\end{displaymath} 

\begin{displaymath}
G_{n,k,2,i}  =  + \frac {b_{3}}{i}\ G_{n-1,k,2,i-1} \quad for \quad i=v_{k,2}+n
\end{displaymath}
\begin{displaymath}
\begin{split} 
G_{n,k,2,i} - t_{p}(i+1)G_{n,k,2,i+1} =  &+ \frac {b_{3}}{i}\
G_{n-1,k,2,i-1}+  b_{4}G_{n-1,k,2,i}\\ &for \quad i=v_{k,2}+n-1
\end{split}
\end{displaymath}\begin{displaymath}
\begin{split} 
G_{n,k,2,i} &- t_{p}(i+1)G_{n,k,2,i+1} =  + \frac {b_{3}}{i}\ G_{n-1,k,2,i-1}+  b_{4}G_{n-1,k,2,i} \\  &+ b_{5}(i+1)G_{n-1,k,2,i+1} \quad for \quad i<=v_{k,2}+n-2
\end{split}
\end{displaymath}

where
\begin{displaymath}
b_{3} = b_{0}- \frac {b_{1}}{t_{p}}\ +  \frac {b_{2}}{t_{p}^{2}}\ 
\end{displaymath}
\begin{displaymath}
b_{4} = b_{1}- \frac {2 \,b_{2}}{t_{p}}\ 
\end{displaymath}
\begin{displaymath}
b_{5} = b_{2}
\end{displaymath}

\item[b)] $i=0 \, k=1$
\begin{displaymath}
\begin{split} 
G_{n,k,1,i} =&+y_{n-1,q}(1)+t_{p} \frac{dy_{n-1,q}(1)}{dt_{n}}\ \\ &- t_{p}G_{n,k,1,i+1} - t_{p}G_{n,k,2,i+1} 
\end{split} 
\end{displaymath}
\begin{displaymath}
\begin{split} 
G_{n,k,2,i} =&-t_{p} \frac{dy_{n-1,q}(1)}{dt_{n}}\ \\ &+ t_{p}G_{n,k,1,i+1} + t_{p}G_{n,k,2,i+1} 
\end{split} 
\end{displaymath}
where $G_{n,k,1,i+1}$ and $G_{n,k,2,i+1}$ exist respectively only for $n+v_{k,1}>=1$ and  $n+v_{k,2}>=1$.

\item[c)] $i=0 \, k>1$
\begin{displaymath}
\begin{split} 
G_{n,k,1,i} =&+y_{n,k-1}( \tau_{k-1})+t_{p}
\frac{dy_{n,k-1}( \tau_{k-1})}{dt_{n}}\ \\ &- \sum_{h=1}^{h=v_{k,1}+n} ( \tau_{k-1}+t_{p} \,h) \tau_{k-1}^{h-1}G_{n,k,1,h}\\ &- e^{- \tau_{k-1}/t_{p}} \sum_{h=1}^{h=v_{k,2}+n}t_{p} \,h \, \tau_{k-1}^{h-1} G_{n,k,2,h}
\end{split} 
\end{displaymath}
\begin{displaymath}
\begin{split} 
G_{n,k,2,i} =&-t_{p}e^{+ \tau_{k-1}/t_{p}} \frac{dy_{n,k-1}( \tau_{k-1})}{dt_{n}}\ \\ &+ e^{+ \tau_{k-1}/t_{p}} \sum_{h=1}^{h=v_{k,1}+n}t_{p} \,h \, \tau_{k-1}^{h-1} G_{n,k,1,h}\\ &+ \sum_{h=1}^{h=v_{k,2}+n} (- \tau_{k-1}+t_{p} \,h) \tau_{k-1}^{h-1}G_{n,k,2,h}
\end{split} 
\end{displaymath}
\end{itemize}

\subsection{Example no.2} 
Since $q=1$, $y_{0,1}=1$, $f_{n,1}=1$ for $n=0$ and $f_{n,1}=0$ for $n>0$,  (\ref{eq:2.4}) and  (\ref{eq:2.1}) become
\begin{equation}\label{eq:3.3}
\begin{split}
\frac{dy_{1,1}(t_{n})} {dt_{n}}\ &+t_{p} \frac{d^{2}y_{1,1}(t_{n})}{dt_{n}^{2}}\ =0 \\ 
\frac{dy_{n,1}(t_{n})} {dt_{n}}\ &+t_{p} \frac{d^{2}y_{n,1}(t_{n})}{dt_{n}^{2}}\ =- b_{0} y_{n-1,1}(t_{n})\\ -b_{1} \frac{dy_{n-1,1}(t_{n})}{dt_{n}}\ &-b_{2} \frac{d^{2}y_{n-1,1}(t_{n})} {dt_{n}^{2}}\ \quad for \, n>1
\end{split}
\end{equation}
\begin{equation}\label{eq:3.4}
\begin{split}
y_{1,1}(t_{1}) &= G_{1,1,1,0} \\
y_{n,1}(t_{n}) &=  \sum_{i=0}^{i=n-1} G_{n,1,1,i}t_{n}^{i}+ e^{-t_{n}/t_{p}} \sum_{i=0}^{i=n-2} G_{n,1,2,i}t_{n}^{i} \quad for \, n>1
\end{split}
\end{equation}
Comparing (\ref{eq:3.1}) and (\ref{eq:3.2}) with (\ref{eq:3.3}) and (\ref{eq:3.4}), it follows that $G_{1,1,1,0}=1$ and the expressions of the coefficients $G_{n,1,p,i}$, detailed in Example no.1 for $i>0$ and $i=0 \, k=1$, are valid also for $n>1$ in Example no. 2 assuming $q=1$, $v_{1,1}=-1$ and $v_{1,2}=-2$. 

\section{Conclusion}
The found explicit solution of neutral differential difference equations will allow to refine the existing  tuning procedures and to implement new ones for the PID or other controllers, used in linear systems of unlimited order provided with one time delay.

\appendix
\section{Derivatives of $y=e^{r \, t}t^{i}$}
The hth derivative of $y=e^{r \, t}t^{i}$ is given by 
\begin{equation}\label{eq:A.1}
\begin{split}
\frac{d^{h}y}{dt^{h}} &= e^{r \, t}  \sum_{j=i-h}^{j=i} \frac {i!h!}{j!(i-j)!(h-i+j)!}\ r^{h-i+j} t^{j} \quad for \, h \leq i\\ \frac{d^{h}y}{dt^{h}} &= e^{r \, t}  \sum_{j=0}^{j=i} \frac {i!h!}{j!(i-j)!(h-i+j)!}\ r^{h-i+j} t^{j} \quad for \, h > i
\end{split}
\end{equation}
These formulas can be checked by comparing the derivative of (\ref{eq:A.1}) to the expression assumed by (\ref{eq:A.1}) when $h$ is replaced with $h+1$.

\section{Identities of series products}
Let define the series products $S_{1}(m)$ and $S_{2}(m,z)$ as
\begin{equation}\label{eq:B.1}
S_{1}(m) = \sum_{h=1}^{h=m} \sum_{i=0}^{i=h-1} \sum_{j=0}^{j=i}
\end{equation}
\begin{equation}\label{eq:B.2}
S_{2}(m,z) = \sum_{h=1}^{h=m} \sum_{i=h}^{i=z} \sum_{j=i-h}^{j=i}
\end{equation}
Shifting in (\ref{eq:B.1}) and (\ref{eq:B.2}) the $j$ series from the third to the second position, we obtain 
\begin{equation}\label{eq:B.3}
S_{1}(m) =  \sum_{h=1}^{h=m} \sum_{j=0}^{j=h-1} \sum_{i=j}^{i=h-1}
\end{equation}
\begin{equation}\label{eq:B.4}
S_{2}(m,z) = S_{21}+S_{22}+S_{23}+S_{24}+S_{25}+S_{26}
\end{equation}
where
\begin{displaymath}
S_{21} = \sum_{h=1}^{h=int[z/2]} \sum_{j=0}^{j=h-1} \sum_{i=h}^{i=h+j}
\end{displaymath}
\begin{displaymath}
S_{22}= \sum_{h=1}^{h=int[z/2]} \sum_{j=h}^{j=z-h} \sum_{i=j}^{i=h+j} 
\end{displaymath}
\begin{displaymath}
S_{23} = \sum_{h=1}^{h=int[z/2]} \sum_{j=z-h+1}^{j=z} \sum_{i=j}^{i=z}
\end{displaymath}
\begin{displaymath}
S_{24} = \sum_{h=int[z/2]+1}^{h=m} \sum_{j=0}^{j=z-h-1} \sum_{i=h}^{i=h+j}
\end{displaymath}
\begin{displaymath}
S_{25}=  \sum_{h=int[z/2]+1}^{h=m} \sum_{j=z-h}^{j=h} \sum_{i=h}^{i=z}
\end{displaymath}
\begin{displaymath}
S_{26} =  \sum_{h=int[z/2]+1}^{h=m} \sum_{j=h+1}^{j=z} \sum_{i=j}^{i=z}
\end{displaymath}
Shifting in (\ref{eq:B.3}) and in (\ref{eq:B.4}) the $j$ series from the second to the first position, we obtain 
\begin{displaymath}
S_{1}(m) =  \sum_{j=0}^{j=m-1}  \sum_{h=j+1}^{h=m} \sum_{i=j}^{i=h-1}
\end{displaymath}
\begin{displaymath}
S_{21} = \sum_{j=0}^{j=int[z/2]-1} \sum_{h=j+1}^{h=int[z/2]} \sum_{i=h}^{i=h+j}
\end{displaymath}
\begin{displaymath}
\begin{split}
S_{22} = &+ \sum_{j=1}^{j=int[z/2]} \sum_{h=1}^{h=j} \sum_{i=j}^{i=h+j} + \sum_{j=int[z/2]+1}^{j=z-1} \sum_{h=1}^{h=z-j} \sum_{i=j}^{i=h+j}\\ 
\end{split}
\end{displaymath}
\begin{displaymath}
S_{23} = \sum_{j=z+1-int[z/2]}^{j=z} \sum_{h=z-j+1}^{h=int[z/2]} \sum_{i=j}^{i=z}
\end{displaymath}
\begin{displaymath}
S_{24} = \sum_{j=0}^{j=z-1-m} \sum_{h=int[z/2]+1}^{h=m} \sum_{i=h}^{i=h+j} + \sum_{j=z-m}^{j=z-2-int[z/2]} \sum_{h=int[z/2]+1}^{h=z-1-j} \sum_{i=h}^{i=h+j}
\end{displaymath}
\begin{displaymath}
\begin{split}
S_{25}=  &+ \sum_{j=z-m}^{j=z-int[z/2]-1} \sum_{h=z-j}^{h=m} \sum_{i=h}^{i=z} + \sum_{j=z-int[z/2]}^{j=int[z/2]} \sum_{h=int[z/2]+1}^{h=m} \sum_{i=h}^{i=z}\\ &+ \sum_{j=int[z/2]+1}^{j=m} \sum_{h=j}^{h=m} \sum_{i=h}^{i=z}
\end{split}
\end{displaymath}
\begin{displaymath}
S_{26} = \sum_{j=int[z/2]+2}^{j=m+1} \sum_{h=int[z/2]+1}^{h=j-1} \sum_{i=j}^{i=z} + \sum_{j=m+2}^{j=z} \sum_{h=int[z/2]+1}^{h=m} \sum_{i=j}^{i=z}
\end{displaymath}

Each series identity can be easily verified by means of the Terms Grid, detailed as example in Table 1. for (\ref{eq:B.1}) and (\ref{eq:B.3}). The symbol x in the grid cell, whose row and column numbers are respectively $i$ and  $j$, means the existence of the term, whose indexes are $i$ and $j$. 

\begin{table}[htbp]
\centering
\caption{Terms Grid of (\ref{eq:B.1}) and (\ref{eq:B.3})} 
\label{Tab1} 
{\begin{tabular}{|c|c|c|c|c|c|c|}
\hline
\multicolumn{2}{|c|}{} & \multicolumn{5}{|c|}{j} \\
\cline{3-7}
\multicolumn{2}{|c|}{} & 0 & 1 & 2 & & h-1 \\
\hline
   & 0 & x & & & & \\
\cline{2-7}
   & 1 & x & x & & & \\
\cline{2-7}
 i & 2 & x & x & x & & \\
\cline{2-7}
   &  & x & x & x & x & \\
\cline{2-7}
 & h-1 & x & x & x & x & x \\
\hline
\end{tabular}}
\end{table}

\section{Vandermonde equations system}
A Vandermonde equations system of mth order is described by
\begin{equation}\label{eq:C.1}
 \sum_{j=1}^{j=m} d_{j} r_{j}^{i-1} x_{i} = Z_{i} \quad \ 1 \leq i \leq m
\end{equation}
where $d_{j}$, $r_{j}$ and $Z_{i}$ are constant.
The solution, evaluated with the Cramer's metod, is given by
\begin{equation}\label{eq:C.2}
x_{j} = \frac { \sum_{i=1}^{i=m} (-1)^{i+j} U_{i,j} Z_{i}}{V}
\end{equation}
where $V$ is the determinant of the coefficient matrix $M_{V}$  and $U_{i,j}$ is the determinant of the submatrix $M_{U}$,  obtained deleting the ith row and the jth column from this matrix. The determinant $V$ is obviously a Vandermonde one and also $U_{i,j}$ can be expressed as a sum of  products of two Vandermonde submatrix determinants, if it is calculated according to the Laplace's expansion, i. e.  assuming
\begin{itemize}
\item[-] for the first submatrix $M_{1}$\\
as rows  the $i-1$ rows from $1$ to $i-1$ of $M_{V}$ and as columns  the $i-1$ columns, belonging to one of the $N=(m-1)!/((i-1)!(m-i)!)$ sets and selected among the $m-1$ ones of  $M_{U}$.
\item[-] for the second submatrix $M_{2}$\\
as rows  the $m-i$ rows from $i+1$ to $m$ of $M_{V}$ and  as columns the remaining ones of $M_{U}$ i.e. not included in $M_{1}$.
\end{itemize}

In detail:
\begin{equation}\label{eq:C.3}
V = \prod_{a>b}(r_{a} -r_{b}) \prod_{j=1}^{j=m} d_{j} 
\end{equation}
\begin{equation}\label{eq:C.4}
U_{i,j} =  \sum_{n=1}^{n=N} (-1)^{s_{k}} W_{i,j,n}
\end{equation}
where
\begin{displaymath}
s_{n} =  \sum_{h=1}^{h=i-1} h+ \sum_{k=1}^{k=i-1} pa(n,k)
\end{displaymath}
\begin{displaymath}
\begin{split}
 W_{i,j,n} =&  \prod_{pa(n,k)>pb(n,k)}(r_{pa(n,k)} -r_{pb(n,k)})  \prod_{k=1}^{k=i-1} d_{pa(n,k)}\\  &\prod_{qa(n,k)>qb(n,k)}(r_{qa(n,k)} -r_{qb(n,k)})  \prod_{k=1}^{k=m-i} d_{qa(n,k)} r_{qa(n,k)}^{i}
\end{split}
\end{displaymath}
\begin{itemize}
\item $h$ is the row number of $M_{V}$
\item  $k$ is the current column number of $M_{1}$ or $M_{2}$
\item $pa(n,k)$ or $pb(n,k)$ is the number of the column, included as kth in $M_{1}$ in increasing order, of $M_{U}$ 
\item $qa(n,k)$ or $qb(n,k)$ is the number of the column, included as kth in $M_{2}$ in increasing order, of $M_{U}$ 
\end{itemize}

\end{document}